\documentclass[final,3p,mathptmx]{elsarticle}
\usepackage{amssymb, amsmath, amsthm}
\usepackage{tikz}
\usepackage{epstopdf}

\biboptions{sort&compress}

\newtheorem{theorem}{Theorem}
\newtheorem{proposition}[theorem]{Proposition}
\newtheorem{lemma}[theorem]{Lemma}
\newtheorem{corollary}[theorem]{Corollary}

\newtheorem{remark}[theorem]{Remark}
\newtheorem{definition}[theorem]{Definition}

\newcommand\bd{\mathbf{d}}

\newcommand\be{\mathbf{e}}

\newcommand\bg{\mathbf{g}}

\newcommand\ba{\mathbf{a}}
\newcommand\bff{\mathbf{f}}
\newcommand\bb{\mathbf{b}}
\newcommand\bh{\mathbf{h}}

\newcommand\ZZ{{\mathbb Z}}
\newcommand{\RR}{\mathbb{R}}
\newcommand{\NN}{\mathbb{N}}
\newcommand{\CC}{\mathbb{C}}

\newcommand{\eop}{$\qquad \Box$}

\def\bdelta{\mbox{\boldmath $\delta$}}

\def\bphi{\mbox{\boldmath $\phi$}}

\numberwithin{equation}{section}

\def\low#1{_{_{{#1}}}}

\newcommand{\baa}{\begin{eqnarray*}}
\newcommand{\eaa}{\end{eqnarray*}}
\newcommand{\ban}{\begin{eqnarray}}
\newcommand{\ean}{\end{eqnarray}}
\newcommand{\PD}{\mathop{\prod}\limits}
\newcommand{\SPAN}{\mathop{\mathrm{span}}\nolimits}

\def\ak{{a}^{[k]}}

\def\aki{{\rm a}_i^{[k]}}
\def\bak{{\bf a}^{[k]}}

\def\NS{S_{\bak}}

\def\vp{\varphi}

\def\mP{\mathcal P}

\newcommand\ra{{\rm{a}}}
\newcommand\rb{{\rm{b}}}

\newcommand\rf{{\rm{f}}}

\newcommand\re{{\rm{e}}}
\newcommand\bbf{\mathbf{f}}

\newcommand{\SUMi}{\sum_{i\in\ZZ}}

\newcommand{\mfrac}[2]%
{\raisebox{0.5pt}{\footnotesize$\dfrac{#1}{#2}$}}
\newcommand{\mbinom}[2]%
{\raisebox{0.5pt}{\footnotesize$\dbinom{#1}{#2}$}}
\def\smmat\{#1&#2\cr#3&#4\}%
   {\!\!\pmatrix{\!#1\!&\!#2\!\cr\!#3\!&\!#4\!}\!\!}
\newcommand\scrm{{\raise0.5pt\hbox{-}}}
\def\eop{{ \vrule height7pt width7pt depth0pt}\par\bigskip} 
\newcommand\ie{{\it\thinspace i.e.}}
{\par\noindent\textbf{Proof:}~}

\begin{document}

\title{Approximation order and approximate sum rules in subdivision}
\author[cc]{Costanza Conti}
\ead{costanza.conti@unifi.it}

\author[lr]{Lucia Romani}
\ead{lucia.romani@unimib.it}

\author[jy]{Jungho Yoon\corref{cor1}}
\ead{yoon@ewha.ac.kr}

\cortext[cor1]{Corresponding author}
\address[cc]{Dipartimento di Ingegneria Industriale, Universit\`{a} di Firenze, Viale Morgagni 40/44, 50134 Italy}
\address[lr]{Dipartimento di Matematica e Applicazioni, Universit\`{a} di Milano-Bicocca, Via R. Cozzi 55, 20125 Milano, Italy  }
\address[jy]{Department of Mathematics, Ewha W. University Seoul, South Korea
}

\begin{abstract}

Several properties of stationary subdivision schemes are nowadays well understood. In particular, it is known that the polynomial generation and reproduction capability of a stationary subdivision scheme is strongly connected with sum rules, its convergence, smoothness and approximation order.
The aim of this paper is to show that, in the non-stationary case, exponential polynomials and approximate sum rules
play an analogous role of polynomials and sum rules in the stationary case.
Indeed, in the non-stationary univariate case we are able to show the following important facts: i) reproduction of $N$ exponential polynomials implies approximate sum rules of order $N$;
ii) generation of $N$ exponential polynomials implies approximate sum rules of order $N$, under the additional assumption of asymptotical similarity and reproduction of one exponential polynomial;
iii) reproduction of an $N$-dimensional space of exponential polynomials and asymptotical similarity imply approximation order $N$; iv) the sequence of basic limit functions of a non-stationary scheme reproducing one exponential polynomial converges uniformly to the basic limit function of the asymptotically similar stationary scheme.
\end{abstract}

\maketitle




\section{Introduction}
In this paper we investigate theoretical properties of non-stationary subdivision schemes.
In particular, we study their approximation order and the role played by approximate sum rules, a non-stationary extension of the well-known notion of sum rules. The obtained results allow us to point out similarities and differences between the stationary and the non-stationary cases.

\smallskip \noindent
For unfamiliar readers, we briefly introduce
subdivision schemes as efficient tools to design smooth curves and surfaces out of sequences of initial points.
As a matter of fact, a subdivision curve or surface is obtained as the limit of an iterative procedure based on the repeated application of local refinement rules generating denser and denser sets of points starting from a coarse initial set roughly describing the desired limit shape \cite{CaravettaDahmenMicchelli,DL2002,Warren2002}.
In practical use, however, only a limited number of iterations are needed.
As a consequence, subdivision schemes are very efficient if compared with traditional parametric curve and surface
representations. They also stand out for ease of implementation and versatility in building free-form surfaces of arbitrary
topological genus.
All these advantages are the reasons for the overwhelming development
of subdivision methods and their increasing use in many applicative areas
such as computer--aided geometric design \cite{DL2002,Warren2002},
curve and surface reconstruction \cite{PR08}, wavelets and multiresolution
analysis \cite{Chui10}, signal/image processing \cite{CRU14,Uhlmann14},
computer games and animation \cite{DeRose2001}.
Within the variety of subdivision methods studied in the literature,
the class of non-stationary subdivision schemes is
currently receiving great attention. This is due to the fact that
non-stationary subdivision schemes are general and flexible enough to
overcome the restricted capabilities of stationary subdivision schemes.
As an example, we can think of the fact that stationary subdivision schemes
are not capable of representing conic sections or, in
general, exponential polynomials.
On the contrary, non-stationary schemes can also generate exponential
polynomials or exponential B-splines,
that is piecewisely defined exponential polynomials
\cite{BCR07,CharinaContiRomani14,ContiRomaniGemignai10,DLL03, JLY, LY10, MWW01,NR15,R10,Schaefer03}.
Reproduction of piecewise exponential polynomials is important in several applications, e.g., in biomedical imaging, in geometric design
and in isogeometric analysis. Moreover, non-stationary subdivision
schemes include Hermite subdivision schemes. Hermite subdivision schemes
are  iterative methods mapping, at each iteration, a set of vector
data consisting of functional values and associated derivatives,
to a denser set of vector data of the same type \cite{MerSau,CMR14}.
They are applied  in geometric modeling for the construction of curves
and surfaces out of points and directional derivatives, and have
recently found application in other contexts such as, for example,
in the design of one-step numerical methods for the numerical solution
of ODE Initial Value Problems.

\smallskip \noindent
The main goal of this paper is to investigate the approximation order of non-stationary subdivision schemes and the role played by approximate sum rules, a non-stationary extension of the well-known notion of sum rules.
Approximate sum rules allow us to link the response of non-stationary  subdivision schemes to specific types of initial data
(precisely data sampled from exponential polynomials) with the approximation and smoothness orders. With some extent we find that exponential polynomials  and approximate sum rules
play an analogous role of polynomials and sum rules in the stationary case. However, important differences unexpectedly arise.

\smallskip \noindent
In stationary subdivision, it is well-known that:
\begin{itemize}
\item The property of \emph{reproducing polynomials of degree $N-1$} (\ie, the capability of a stationary subdivision scheme to reproduce in the limit exactly the same degree-$(N-1)$ polynomial from which the data are sampled) implies sum rules of order $N$ (see \cite{CharinaConti12}). Moreover, it implies \emph{approximation order $N$} (see \cite{ALevin2003}).
\item The property of \emph{generating polynomials of degree $N-1$} (\ie, the capability of a stationary subdivision scheme to provide polynomials of degree $N-1$ as limit functions) is equivalent to the fulfillment of sum rules of order $N$ (see \cite{CharinaConti12}). Moreover, it guarantees the existence of difference operators whose spectral properties characterize convergence and regularity of the subdivision scheme (see \cite{DL2002}).
\item \emph{Sum rules of order $N$} are a necessary condition for convergence and for $C^{N-1}$-continuity of stationary subdivision schemes (see, e.g., \cite{Cabrelli,CharinaContiSauer2005,Han2008,Han2010,JiaHan2007} and references therein).
\end{itemize}

\smallskip \noindent
In contrast, in the non-stationary setting, although the notions of generation and reproduction of polynomials are straightforwardly replaced by the notions of generation and reproduction of exponential polynomials, the situation is not so clear.
For this reason, in this paper, we prove the following important results, which allow us to point out similarities and differences between the stationary and non-stationary cases:

\smallskip \noindent
\begin{itemize}
\item As in the stationary case, the property of \emph{reproducing $N$ exponential polynomials} implies approximate sum rules of order $N$.
Moreover, we are able to show that it implies \emph{approximation order $N$} if asymptotical similarity to a convergent stationary scheme is assumed.
\item The property of \emph{generating $N$ exponential polynomials} implies \emph{approximate sum rules of order $N$} if asymptotical similarity to a convergent stationary scheme and reproduction of one exponential polynomial are assumed.
Moreover, as in the stationary case, the property of generating exponential polynomials guarantees the existence of \emph{difference operators}.
\item \emph{Approximate sum rules of order $N$} and asymptotical similarity to a stationary $C^{N-1}$ subdivision scheme provide sufficient conditions for $C^{N-1}$ regularity of non-stationary subdivision schemes.
\end{itemize}

\smallskip \noindent
We emphasize that, in order to maintain the notation and the proofs of the results as simple as possible, we deliberately focus the discussion on univariate non-stationary subdivision schemes only.
Note also that, motivated by the fact that primal and dual subdivision schemes are essentially the ones of interest in applications, we restrict our attention to these two subclasses of subdivision schemes, which are known to include odd and even symmetric subdivision schemes.

\medskip \noindent
The remainder of this paper is  organized as follows.
In Section \ref{sec:basic}, after providing some basic definitions, we recall known results concerning exponential polynomial generation and reproduction. Then, in Section \ref{sec:link}
we study the link between exponential polynomial generation/reproduction and approximate sum rules.
Section \ref{sec:asymptotic_blf} focuses on the uniform convergence of the sequence of basic limit functions of a non-stationary scheme
to the basic limit function of the asymptotically similar stationary scheme. This extends an existing result in \cite{DynLevin95}, which resorts to the assumption of asymptotical equivalence between the two schemes. Finally,
under the weaker assumption of asymptotical similarity rather than asymptotical equivalence to a stationary scheme,
in Section \ref{Sec:approx} we also investigate the approximation order of a non-stationary subdivision scheme reproducing a
space of exponential polynomials.
All the results of  this paper are summarized in Section \ref{conclusions}.

\section{Symbols, exponential polynomial generation and reproduction}\label{sec:basic}

Given an initial set of discrete data $\bff^{[0]}=\{ \rf_i^{[0]},\ i\in \ZZ\}$, a univariate, non-stationary subdivision scheme constructs
the sequence of refined data  $\{\bff^{[k]},\ k> 0\}$ via the repeated application of the subdivision operators $S_{\ba^{[0]}},\cdots,S_{\ba^{[k-1]}}$
associated with the finitely supported sequences of coefficients $\{\ba^{[k-1]},\ k > 0\}$, $\ba^{[k]}:= \{\ra_i^{[k]} \in \RR,\ i\in \ZZ\}$, named the \emph{subdivision masks}. We assume that all masks have the same support $\{-M,\cdots,M\}$, $M\in\NN$, \ie, for all $k \geq 0$, $\ra_i^{[k]}=0$ if $i<-M$ or $i>M$. To generate the refined data sequence $\bff^{[k]}=\{\rf^{[k]}_i, \, i \in \ZZ\}$, $k>0$, the subdivision operator $S_{\ba^{[k-1]}}$ is applied to the
data sequence $\bff^{[k-1]}=S_{\ba^{[k-2]}} \cdots S_{\ba^{[0]}} \bff^{[0]}$ acting as
\begin{equation}\label{NS-SUB}
\rf^{[k]}_i:=\left(S_{\ba^{[k-1]}}\bff^{[k-1]}\right)_i
:=\sum_{j\in \ZZ}\ra^{[k-1]}_{i-2j} \, \rf_j^{[k-1]}, \quad i\in \ZZ.
\end{equation}
A non-stationary subdivision scheme is thus identified by the collection of subdivision operators
and therefore denoted as $\{S_{\ba^{[k]}},\ k \geq 0\}$.
When all subdivision operators are the same, the subdivision scheme is called
\emph{stationary} and is simply denoted as $\{S_{\ba}\}$.

\begin{definition}
A subdivision scheme is said to be $C^\ell$-convergent
if, for any initial sequence $\bff^{[0]}\in \ell_\infty(\ZZ)$, there exists
a function $g_{\bff^{[0]}}\in C^\ell(\RR)$ such that
$$
\lim_{k \rightarrow \infty}
\sup_{i \in \ZZ} \vert g_{\bff^{[0]}}(i2^{-k})- (S_{\ba^{[k-1]}} \cdots S_{\ba^{[0]}}\bff^{[0]})_i   \vert =0\,.
$$
If $\ell = 0$, the scheme is simply said to be convergent.
\end{definition}

\smallskip \noindent
For the practical use of a convergent subdivision scheme it is important to link the properties of the limit function $g_{\bff^{[0]}}$ with the properties of the initial sequence $\bff^{[0]}$. This is particularly true in the case when $\bff^{[0]}$ consists of samples of special types of functions: polynomial, trigonometric and hyperbolic functions or, more generally, exponential polynomial functions.
Indeed, the response of the subdivision scheme to these types of starting data is not only important for the design of shapes
of practical interest in applications (see, e.g., \cite{BCR07,R09,Schaefer03,Warren2002}), but is also strongly connected to the following key properties of the subdivision scheme: its smoothness and its approximation order (see \cite{Cabrelli,CharinaContiSauer2005,CharinaContiGuglielmiProtasov2014,ChenJia,Han2008,Han2010,JetterPlonka,JiaHan2007,Jia95,JiaJiang2002} and Section \ref{sec:link} as well as \cite{JKLY} and Section \ref{Sec:approx}, respectively).

\smallskip \noindent
Throughout this paper we use the term \emph{generation} to refer to the capability of a subdivision scheme of providing specific types of limit functions. On the contrary, by the term \emph{reproduction} we refer to the capability of a subdivision scheme of reproducing in the limit exactly the same function from which the data are sampled.\\
To study the capabilities of a non-stationary subdivision scheme $\{S_{\ba^{[k]}}, \ k \geq 0\}$ of generating/reproducing
exponential polynomial functions, we need the following definition.

\begin{definition}\label{def:space-expol}
Let $x\in \RR$ and $\eta,N \in \NN$.
Assume also $\lambda_n\in \CC$ $n=1,\cdots,\eta$, and $\mu_n \in \NN$ $n=1,\cdots,\eta$,  to be such that $\sum_{n=1}^{\eta}\mu_n=N$. The $N$-dimensional space    spanned by
$N$ linearly independent exponential polynomials is $$\Phi_N:={\rm span}\{x^{\beta} e^{\lambda_n x}, \, \beta=0,\cdots,\mu_n-1, \, n=1,\cdots,\eta\}.$$
\end{definition}

\smallskip \noindent
For simplicity, in the following we write  $\Phi_N$ as $\Phi_N={\rm span}\{\varphi_0(x),\cdots,\varphi_{N-1}(x)\}$.\\

The next definition stresses the difference between the notions of generation and reproduction of the space of exponential polynomials $\Phi_N$.

\begin{definition}\label{def:ERgenerationlimit}
Let $t^{[0]}_i$, $i\in \ZZ$, be ordered values on the real axis such that $t^{[0]}_{i+1}-t^{[0]}_i=1$ for all $i \in \ZZ$.
A convergent, non-stationary subdivision scheme $\{S_{\ba^{[k]}}, \ k \geq 0\}$ is said to be
\begin{itemize}
\item[i)] \emph{$\Phi_N$-generating}, if for all initial sequences ${\bbf}^{[0]}:=\{f(t^{[0]}_i),\ i\in \ZZ\}$, $f\in \Phi_N$, we get
$$\displaystyle{\lim_{k\rightarrow \infty}S_{\ba^{[k]}} S_{\ba^{[k-1]}} \cdots S_{\ba^{[0]}}{\bbf}^{[0]} \in \Phi_N};$$
\item[ii)] \emph{$\Phi_N$-reproducing}, if for all initial sequences ${\bbf}^{[0]}:=\{f(t^{[0]}_i),\ i\in \ZZ\}$, $f\in \Phi_N$, we get
$$\displaystyle{\lim_{k\rightarrow \infty}S_{\ba^{[k]}} S_{\ba^{[k-1]}} \cdots S_{\ba^{[0]}}{\bbf}^{[0]}=f}.$$
\end{itemize}
\end{definition}

\begin{remark}
It is easy to see that the space of polynomials of degree at most $N-1$, can be obtained from Definition \ref{def:space-expol}
by choosing $\lambda_1=0$ and $\beta=0, \cdots, N-1$.
As a consequence, Definition \ref{def:ERgenerationlimit} applies  also in the stationary case, namely when $\{S_{\ba^{[k]}}, \ k \geq 0\}$ is simply $\{S_{\ba}\}$ and $\Phi_N$ reduces to $\Pi_{N-1}:={\rm span}\{1,x,\cdots,x^{N-1}\}$.
\end{remark}

\smallskip \noindent
For a given non-stationary subdivision scheme $\{S_{\ba^{[k]}}, \ k \geq 0\}$, its generation and reproduction properties
are encoded in the so-called subdivision symbols
$$a^{[k]}(z):=\sum_{i\in \ZZ}\ra^{[k]}_{i}z^i, \quad z \in \CC \setminus \{0\}, \quad k\ge 0,$$
associated to the subdivision masks
$$\ba^{[k]}:=\{\ra^{[k]}_{i} \in \RR, \, i \in \ZZ \}, \quad k \geq 0.$$
Then \eqref{NS-SUB} can be written
as
\begin{equation}\label{F-NS}
f^{[k+1]} (z) = a^{[k]}(z) f^{[k]} (z^2),\quad k\ge 0
\end{equation}
where
$$f^{[k]} (z) =\sum_{i\in \ZZ} \rf^{[k]}_i z^i.$$

\smallskip \noindent
To simplify the presentation of the algebraic conditions on non-stationary  subdivision symbols that guarantee certain generation and reproduction properties of the corresponding subdivision scheme, we start from reviewing the results of the stationary case.\\
It is well-known that a stationary subdivision scheme $\{S_{\ba}\}$ is $\Pi_{N-1}$-generating if and only if its symbol $a(z)$ can be written in the form
\begin{equation}\label{def:rep2}
a(z)=(1+z)^N \, b(z),
\end{equation}
for some Laurent polynomial $b(z)$ (see \cite[Theorem 3.2]{DHSS08}).
In  \cite[Theorems 4.6, 4.7]{DHSS08}  the authors also showed that for a $\Pi_{N-1}$-generating subdivision scheme $\{S_{\ba}\}$ with symbol $a(z)$, a necessary and sufficient condition for reproducing $\Pi_{N-1}$ is given by the existence of
a Laurent polynomial $c(z)$ such that
\begin{equation}\label{def:rep}
2-a(z^{1+\nu}) \, z^{\nu}=(1-z)^N \, c(z) \qquad \hbox{with} \quad \nu \in \{0,1\}.
\end{equation}\\
To better understand the role of the parameter $\nu \in \{0,1\}$ in \eqref{def:rep}, we need to recall the concept of \emph{parametrization} of a subdivision scheme, i.e. the choice of the grid points $t_i^{[k]}$, $i\in \ZZ$, to which the $k$-level data $\rf_i^{[k]}$, $i\in \ZZ$, are attached.
The general definition of the grid points $t_i^{[k]}$, $i\in \ZZ$, is
\begin{equation}\label{eq:param}
t_i^{[k]}:=\frac{i+p}{2^k}, \qquad p:=-\left(\tau+\frac{\nu}{2}\right), \quad \tau \in \ZZ, \quad \nu\in\{0,1\},
\end{equation}
where the number $p$ is called the \emph{shift parameter} (see \cite{CH11}).
When $\nu=0$, then $p \in \ZZ$ and the sets of parameters ${\mathbf t}^{[k]}=\{t_i^{[k]}, i \in \ZZ\}$, $k \geq 0$, provide the so-called \emph{primal} parametrization. Subdivision schemes based on this choice are thus called \emph{primal} subdivision schemes.
In contrast, when $\nu=1$, then $p \in \frac{1}{2}\ZZ$ and the parametrization is called \emph{dual} as well as the corresponding subdivision schemes.
To simplify the notation, it is convenient to assume $p=-\frac{\nu}{2}$, $\nu\in\{0,1\}$, namely $p \in \{0, - \frac{1}{2}\}$.
In fact, if a subdivision scheme is $\Pi_{N-1}$-reproducing with respect to the shift parameter $p$,
we can always multiply its symbol by $z^{\tau}$  to make it $\Pi_{N-1}$-reproducing with respect to the shift parameter
$p+\tau$, namely to $-\frac{\nu}{2}$.

\smallskip \noindent
In  the non-stationary case, conditions on the subdivision symbols that guarantee exponential polynomial generation were first given in \cite{VBU}. In case of  primal and dual schemes only, a complete set of algebraic conditions to be satisfied by the subdivision symbols of a non-stationary scheme in order to guarantee exponential polynomial generation and reproduction, can be found both
in \cite{ContiRomani11} and in \cite{JKLY}. In the following theorem we recall the set of conditions given in the second reference.
They are the non-stationary counterparts of \eqref{def:rep2} and \eqref{def:rep}, and are stated in \eqref{AK-10} and \eqref{EXP-REP}, respectively.

\begin{theorem}\label{THM-E-GEN}
Let $\Phi_N$ denote the $N$-dimensional space of exponential polynomials given in Definition \ref{def:space-expol}.
If the subdivision scheme $\{ \NS,\ k \geq 0 \}$ is \emph{$\Phi_N$-generating},
then its symbols $\{ \ak(z),\ k \geq 0 \}$ satisfy
\begin{eqnarray}\label{AK-10}
 a^{[k]}(z)
=b^{[k]}(z)\PD_{n=1}^{\eta}(1+e^{\lambda_n 2^{-k-1}} z)^{\mu_n}, \quad k\geq 0,
\end{eqnarray}
for some Laurent polynomials $\{ b^{[k]}(z),\ k \geq 0 \}$.
Moreover, if $\{ \NS,\ k \geq 0 \}$ is \emph{$\Phi_N$-reproducing} with respect to the parametrization in \eqref{eq:param}, then
\ban \label{EXP-REP}
2 - a^{[k]}(z^{1+\nu})z^{\nu}
=c^{[k]}(z)\PD_{n=1}^{\eta} (1-e^{\lambda_n 2^{-k-1-\nu}} z)^{\mu_n}, \quad \nu \in \{0,1\},\qquad k\ge 0,
\ean
for some Laurent polynomials $\{ c^{[k]}(z),\ k \geq 0 \}$.
\end{theorem}

\smallskip \noindent
In the stationary case, conditions \eqref{def:rep2} and \eqref{def:rep} simply express the fact that the symbols $a(z)$ and $2 - a(z^{1+\nu})z^{\nu}$ contain the factors $(1+z)^N$ and $(1-z)^N$, respectively.
The analogous conditions for the non-stationary case, stated in equations \eqref{AK-10} and \eqref{EXP-REP}, are nothing but the requirement that the $k$-level symbols $a^{[k]}(z)$ and $2 - a^{[k]}(z^{1+\nu})z^{\nu}$ contain the factors $\PD_{n=1}^{\eta}(1+e^{\lambda_n 2^{-k-1}} z)^{\mu_n}$ and  $\PD_{n=1}^{\eta} (1-e^{\lambda_n 2^{-k-1-\nu}} z)^{\mu_n}$, respectively.

\begin{remark}\label{rem:p}
Condition \eqref{EXP-REP} was proven in \cite{JKLY} under the assumption that $a^{[k]}(z)$ is a symmetric symbol.
However, \eqref{EXP-REP} is valid also for primal/dual subdivision schemes that are not symmetric. In fact, the proof
in \cite{JKLY} actually does not make use of the symmetry assumption.
As an example, consider the not symmetric primal subdivision scheme with $k$-level symbol
$$
a^{[k]}(z)=(1+z)^2 (z+r_k) (z+r_k^{-1}) \big( 1+b_1^{[k]} z+b_2^{[k]} z^2+b_3^{[k]} z^3 \big),
$$
where $r_k= e^{\lambda 2^{-k-1}}$ and
$$
\begin{array}{l}
b_1^{[k]}=-\frac{4 r_k^4 + 6 r_k^3 + 9 r_k^2 + 6 r_k + 4}{2 (r_k^4 + 2 r_k^3 + 2 r_k^2 + 2 r_k + 1)}, \quad
b_2^{[k]}=\frac{2 r_k^2 + r_k + 2}{2(r_k^2 + 1)}, \quad
b_3^{[k]}=-\frac{r_k (2 r_k^2 + r_k + 2)}{2 (r_k^2 + 1) (r_k+1)^2}.
\end{array}
$$
It is not difficult to show that it is convergent and reproduces $\Phi_4=\SPAN\{1,x,e^{\pm \lambda x}\}$ with respect to the primal parametrization $\{t_i^{[k]}=\frac{i}{2^k},\ i\in \ZZ\}$, thus it indeed satisfies \eqref{AK-10} and \eqref{EXP-REP} with $\nu=0$.\\
Similarly, we can also easily see that the not symmetric dual subdivision scheme with $k$-level symbol
$$
a^{[k]}(z)=(1+z) (z+r_k) (z+r_k^{-1}) \big(b_1^{[k]} +b_2^{[k]} z +b_3^{[k]}z^2 \big),
$$
where $r_k= e^{\lambda 2^{-k-1}}$ and
$$
\begin{array}{l}
b_1^{[k]}=\frac{r_k^{13/2} - r_k^4 - r_k^2 + r_k^{-1/2}}{(r_k^2 - 1)^2 (r_k^2 + 1)}, \quad
b_2^{[k]}=-\frac{ r_k^{13/2} + r_k^{11/2} -r_k^5 -2 r_k^3 -r_k + r_k^{1/2} + r_k^{-1/2}}{(r_k^2 - 1)^2 (r_k^2 + 1)}, \quad
b_3^{[k]}=\frac{r_k (r_k^{9/2} -r_k^3 -r_k +r_k^{-1/2})}{(r_k^2 - 1)^2 (r_k^2 + 1)},
\end{array}
$$
reproduces $\Phi_3=\SPAN\{1,e^{\pm \lambda x}\}$ with respect to the dual parametrization $\{t_i^{[k]}=\frac{i-\frac{1}{2}}{2^k},\ i\in \ZZ\}$ and indeed satisfies \eqref{AK-10} and \eqref{EXP-REP} with $\nu=1$.
\end{remark}


\section{Exponential polynomials and approximate sum rules}\label{sec:link}

The notion of \emph{approximate sum rules} was first introduced in \cite{CharinaContiGuglielmiProtasov2014} as a possible generalization of the notion of \emph{sum rules}
to the non-stationary setting. For the reader's convenience both definitions are recalled here.

\begin{definition}\label{def:sum_rules}
Let $N \in \NN$. The sequence of symbols $\{a^{[k]}(z),\ \ k \geq 0\}$ is said to satisfy \emph{sum rules of order $N$} if for all $k \geq 0$ it is verified that
\begin{equation}\label{def:sumrules0}
a^{[k]}(1)=2 \qquad \hbox{and} \qquad \frac{d^\beta a^{[k]}(-1)}{dz^\beta}=0, \ \ \beta=0,\ldots, N-1.
\end{equation}
\end{definition}

\begin{definition}\label{def:approx_sum_rules} Let $N \in \NN$. The sequence of symbols $\{a^{[k]}(z),\ \ k \geq 0\}$ is said to satisfy \emph{approximate sum rules of order $N$} if it fulfills
\begin{equation}\label{deltak1}
\sum_{k=0}^{\infty} |a^{[k]}(1)-2| <\infty
\qquad \hbox{and} \qquad
\sum_{k=0}^{\infty} 2^{\, k(N-1)}\sigma_k < \infty,\quad \hbox{for}\quad \sigma_k\ :=\ \max_{\beta=0,\ldots, N-1} {2^{-k\,\beta}
   \left |\frac {d^\beta {a}^{[k]}(-1)}{dz^\beta} \right|}\,.
\end{equation}
\end{definition}

\smallskip \noindent
In stationary subdivision, it is well-known that generation of $\Pi_{N-1}$ is equivalent to the fact that the subdivision symbol satisfies sum rules of order $N$, while reproduction of $\Pi_{N-1}$ just implies sum rules of order $N$ (see, e.g., \cite{Cabrelli, ChenJia, ContiRomaniGemignai09, Han2010, JetterPlonka}).\\

\smallskip \noindent
The goal of this section is to study the link between generation/reproduction of exponential polynomials and approximate sum rules.
In particular, in Theorem \ref{Y-THM-10} we show that reproduction of the exponential polynomial space $\Phi_N$ implies approximate sum rules of order $N$, if the $N \times N$ Wronskian matrix of $\Phi_N$ defined by
\begin{equation}\label{def:wronskymat}
{\cal W}_{\Phi_N} (x) := \Big(\frac{1}{\beta!} \frac{d^\beta \varphi_{\alpha}(x)}{dx^\beta} , \
                    \alpha, \beta =0,\ldots,N-1 \Big),
\end{equation}
is invertible for all $x$ in a neighborhood of zero.
Moreover, in Theorem \ref{teo:new}, we prove that generation of $\Phi_N$ and reproduction of just one exponential polynomial in $\Phi_N$ imply approximate sum rules of order $N$ if, besides the invertibility of
${\cal W}_{\Phi_N}(x)$ for all $x$ in a neighborhood of zero, we additionally assume asymptotical similarity of the non-stationary scheme to some convergent stationary scheme.\\
Thus, for the sake of completeness, before formulating Theorem \ref{teo:new} we provide the definition of the notion of asymptotical similarity, originally introduced in \cite{ContiDynManniMazure14}
to weaken the notion of asymptotical equivalence proposed earlier in \cite{DynLevin95}.

\begin{definition}\label{def:asympt_simil}
A non-stationary subdivision scheme $\{S_{\ba^{[k]}}, \, k \geq 0\}$ and a stationary one $\{S_{\ba}\}$ are said to be
\emph{asymptotically similar}, respectively  \emph{asymptotically equivalent}, if
the associated sequence of subdivision masks $\{\ba^{[k]},\ k \geq 0\}$ and the subdivision mask $\{\ba\}$ satisfy
$$\lim_{k\rightarrow\infty}\|\ba^{[k]}-\ba\|=0,\quad \hbox{respectively}\quad \sum_{k=0}^\infty\|\ba^{[k]}-\ba\| <\infty\,.$$
\end{definition}
\noindent
Note that, here and in the sequel, $\| \cdot\|$ stands for the infinity norm
of subdivision operators, sequences or functions,
\ie, $\|S_{\ba^{[k]}}\|:=\max\left\{\sum_{i\in \ZZ}|\ra^{[k]}_{2i}|,\
\sum_{i\in \ZZ}|\ra^{[k]}_{2i+1}|\right\}$,
 $\| \bff\|:=\sup_{i\in \ZZ}|\rf_i|$ and $\| F\|:=\sup_{x\in \RR}|F(x)|$.

\begin{theorem}\label{Y-THM-10}
Let $\Phi_N$ denote the $N$-dimensional space of exponential polynomials given in Definition \ref{def:space-expol} and
let $\{\ak(z), k \geq 0 \}$ be the symbols of a $\Phi_N$-reproducing non-stationary subdivision scheme $\{\NS,\ k\geq 0\}$.
If the Wronskian matrix of $\Phi_N$, ${\cal W}_{\Phi_N}(x)$, is invertible for all $x$ in a neighborhood of zero, then for $\beta=0,\cdots, N-1$ it is verified that
\ban\label{SUM-R0}
| \ak(1) -2|=O(2^{-kN}) \qquad \hbox{and} \qquad
\left | \frac {d^\beta \ak(-1)}{dz^\beta} \right |
      =O(2^{-k(N-\beta)}), \qquad k\to \infty.
\ean
\end{theorem}
\proof
We first note that the $\beta$-th derivative of $\ak(z)$
evaluated at $z=-1$ can be expressed as
\ban\label{D-AK}
\frac {d^\beta \ak(-1)}{dz^\beta}
=\sum_{\ell=0}^\beta \theta_{\beta,\ell}\SUMi \aki i^\ell (-1)^i,
 \ean
for some constants $\theta_{\beta,\ell}$, $\ell=0,\cdots,\beta$ such that $\theta_{\beta,0}=\delta_{\beta,0}$.
Thus, to verify  (\ref{SUM-R0}),
it suffices to show that
\begin{equation}\label{SUM-R0_new}
\left |\SUMi \aki-2 \right |=O(2^{-kN}) \qquad \hbox{and} \qquad
\left | \SUMi (-1)^i i^\beta \aki\right | =O(2^{-k(N-\beta)}),\qquad k\to \infty.
\end{equation}
Our approach for this task is to separate the summation
in equations \eqref{SUM-R0_new} into two parts, and then estimate them separately.
To do this, for each $\beta=0,\cdots, N-1$, we define the pair of functions
\ban\label{PHI-20}
{\mP}_{j,\beta} (x) :=\sum_{n=0}^{N-1} m^{[k]}_{j,\beta,n} \, \vp_n (x),\quad  j=0,1,
\ean
where the coefficient vector ${\bf m}^{[k]}_{j, \beta}=(m^{[k]}_{j, \beta,n},\  n=0,\cdots, N-1)$ is
obtained by solving the Hermite interpolation problem
\ban\label{PHI-0}
\frac{d^\ell \mP_{j,\beta}(j 2^{-k-1})}{dx^\ell}=
 \delta_{\beta,\ell}(-1)^{ \ell} \ell !,\quad \ell= 0,\ldots, N-1,\quad j=0,1.
\ean
Denoting by $\delta_{\beta,\ell}$ the Kronecker symbol, equations \eqref{PHI-0}
translate into the linear systems
\ban\label{WR-10}
{\cal W}_{\Phi_N} (j 2^{-k-1}) \cdot ({\bf m}^{[k]}_{j,\beta})^T
 =  {\bf c}_\beta^T, \quad j=0,1
\ean
with ${\bf c}_\beta :=(\delta_{\beta,\ell }(-1)^{ \ell}\ell !,\ \ell=0,\cdots,N-1)$. Due to the assumption that
${\cal W}_{\Phi_N}(x)$ is invertible for all $x$ in a neighborhood of zero, such a linear system clearly has a unique solution.
Next we define \ban\label{ANU}
 Q_{j, \beta}
   := \sum_{i\in \ZZ}\ra^{[k]}_{j-2i}\
   \Big ((j 2^{-1}- i) 2^{-k}\Big )^\beta - \mP_{j,\beta} (j 2^{-k-1}),
\quad j=0,1,
\ean
and note that, in view of \eqref{PHI-0},
\begin{equation}\label{eq:p0bp1b}
\mP_{0,\beta} (0) = \mP_{1,\beta} (2^{-k-1})=\delta_{\beta,0}.
\end{equation}
Thus, noting that $\mP_{0,0} (0) = \mP_{1,0} (2^{-k-1})=1$,
we get
\begin{equation}\label{ASR-B0}
\SUMi \aki-2
   = Q_{0,0} + Q_{1,0}.
\end{equation}
Moreover, due to \eqref{eq:p0bp1b}, we have
\begin{align}\label{ASR-B}
\begin{split}
2^{-\beta (k+1)} \SUMi (-1)^i i^\beta \aki
&= \SUMi \ra^{[k]}_{-2i} (-i2^{-k})^\beta
   - \SUMi \ra^{[k]}_{1-2i} \big ((2^{-1} -i)2^{-k}\big)^\beta
 =  Q_{0,\beta} - Q_{1,\beta}.
\end{split}
\end{align}
We first estimate $Q_{0,\beta}$. The other term $Q_{1,\beta}$
can be handled analogously.
Since $ \mP_{0 ,\beta}$ is  a linear combination of exponential polynomials in
$\Phi_N$ and the non-stationary scheme $\{\NS,\ k\geq 0\}$ reproduces $\Phi_N$,
we get the identity $\mP_{0,\beta} (0)
= \sum_{i\in \ZZ}\ra^{[k]}_{-2i}\, \mP_{0,\beta}(i2^{-k})$.
Plugging this into (\ref{ANU}) for $j=0$ leads to
\ban\label{QB-10}
Q_{0,\beta}
  =\sum_{i\in \ZZ } \ra^{[k]}_{-2i}
     \Big ( (-i2^{-k})^\beta - \mP_{0,\beta}(i2^{-k}) \Big ).
\ean
Here, we will use the arguments of Taylor expansion for $\mP_{0,\beta}$.
Precisely, let  $T_{\mP_{0,\beta}}^N$ be the degree-$(N-1)$ Taylor
polynomial of the function $\mP_{0,\beta}$
around $0$, that is,
$T_{\mP_{0,\beta}}^N(x) :=
\sum_{\ell=0}^{N-1}\frac{x ^\ell}{\ell !}\,
 \frac{d^\ell \mP_{0,\beta}(0)}{dx^\ell}$.
Then, we replace $\mP_{0,\beta}$ in $Q_{0,\beta}$ by its Taylor
polynomial $T_{\mP_{0,\beta}}^N$ plus the remainder term, say $R_{\mP_{0,\beta}}^N$,
such that we have the form
$
\mP_{0,\beta} (i2^{-k})
=T_{\mP_{0,\beta}}^N(i2^{-k}) + R_{\mP_{0,\beta}}^N(i2^{-k}).
$
In fact, from the Hermite interpolation conditions in (\ref{PHI-0}),
we find that
$T_{\mP_{0,\beta}}^N (i2^{-k})=(-i2^{-k})^\beta.$
Hence, it leads to the equations
\ban\label{QB-20}
Q_{0,\beta}  = - \sum_{i\in \ZZ} \ra^{[k]}_{-2i}
R_{\mP_{0,\beta}}^N(i2^{-k}).
\ean
On the other hand, by \eqref{WR-10},
$ ({\bf m}^{[k]}_{0,\beta})^T
= {\cal W}_{\Phi_N} (0)^{-1} {\bf c}_\beta^T$
such that for a given $\beta$,
the coefficient vector
${\bf m}^{[k]}_{0,\beta}$ for $\mP_{0,\beta}$ in (\ref{PHI-20})
can be bounded independently of $k \geq 0$.
It implies that
the $\alpha$-th derivative of $\mP_{0,\beta}$
for each $\alpha=0,\ldots, N-1$,
is uniformly bounded around the origin.
Consequently,  we get $|R_{\mP_{0,\beta}}^N(i2^{-k})| = O(2^{-kN})$
and hence, by \eqref{QB-20},
$|Q_{0,\beta}| = O(2^{-kN})$ as $k\to \infty$.
Similarly, we can prove  the same convergence rate for $Q_{1,\beta}$, namely
$|Q_{1,\beta}| = O(2^{-kN})$ as $k\to \infty$.
Combining these two convergence properties and applying equations
(\ref{ASR-B0}) and (\ref{ASR-B}), we finally get \eqref{SUM-R0_new}.
Thus, referring back to the identity in (\ref{D-AK}), the proof is completed.
\endproof

\begin{corollary}
Let $\{\NS,\ k\geq 0\}$ be a $\Phi_N$-reproducing non-stationary subdivision scheme with  ${\cal W}_{\Phi_N}(x)$ invertible for all $x$ in a neighborhood of zero. Then $\{\NS,\ k\geq 0\}$ satisfies approximate sum rules of order $N$.
\end{corollary}

\begin{remark}\label{rem:ASR_notsuff}
It is easy to see that conditions \eqref{SUM-R0} are not sufficient for
the reproduction of any exponential polynomial. As a counterexample, consider the level-dependent
perturbation of quadratic B-splines (still a dual scheme with $\nu=1$) given by the $k$-level symbol
\begin{equation}
a^{[k]}(z)=\frac14+2^{-k}+\left(\frac34-2^{-k}\right)z+\left(\frac34-2^{-k}\right)z^2+\left(\frac14-2^{-k}\right)z^3,\\
\end{equation}
which satisfies $a^{[k]}(1)=2(1 - 2^{-k})$ and $a^{[k]}(-1)=2^{-k+1}$.

\smallskip \noindent Although \eqref{SUM-R0} are satisfied with $N=1$, no exponential polynomials can be reproduced by this scheme. To see it, it suffices to observe that, for $k=1$ and $\lambda\in \CC$, the complex function $2-a^{[k]}((e^{-\lambda/2^{k+2}})^2)e^{-\lambda/2^{k+2}}$ is always different from $0$.
\end{remark}

\medskip \noindent
In the next theorem we replace the assumption that $\{\NS,\ k\geq 0\}$ is $\Phi_N$-reproducing
with the weaker assumption that it is $\Phi_N$-generating and reproduces only one exponential polynomial in $\Phi_N$;
to compensate the latter weaker conditions, we additionally require asymptotical similarity of $\{\NS,\ k\geq 0\}$ to a convergent stationary scheme $\{S_\ba\}$.
Therefore, we consider subdivision schemes $\{\NS,\ k\geq 0\}$
satisfying the first condition in Definition \ref{def:asympt_simil}.

\begin{theorem}\label{teo:new}
Let $\Phi_N$ denote the $N$-dimensional space of exponential polynomials given in Definition \ref{def:space-expol} and
let $\{\ak(z), k\geq 0\}$ be the symbols
of a non-stationary subdivision scheme $\{\NS,\ k\geq 0\}$ which is $\Phi_N$-generating and reproduces one exponential polynomial in $\Phi_N$.
If $\{\NS,\ k\geq 0\}$ is asymptotically similar to a convergent stationary subdivision scheme $\{S_\ba\}$, and the Wronskian matrix ${\cal W}_{\Phi_N}(x)$, defined in (\ref{def:wronskymat}), is invertible for all $x$ in a neighborhood of zero, then for $\beta=0,\cdots, N-1$
\ban\label{SUM-R}
| \ak(1) -2 | =O(2^{-k}) \qquad \hbox{and} \qquad
\left | \frac {d^\beta \ak(-1)}{dz^\beta} \right |
 =O(2^{-k(N-\beta)}), \qquad k\to \infty.
\ean
\end{theorem}

\proof Since $\{\NS,\ k\geq 0\}$ is $\Phi_N$-generating, in view of Theorem \ref{THM-E-GEN} we know that equation (\ref{AK-10}) holds true for all $k \geq 0$.
By the asymptotical similarity assumption, as $k\rightarrow \infty$ the left hand side of (\ref{AK-10}) tends to $a(z)$ (the symbol of the stationary subdivision  scheme $\{S_\ba\}$) and  the term
$\PD_{n=1}^{\eta}(1+e^{\lambda_n 2^{-k-1}} z)^{\mu_n}$, appearing in the right hand side, tends to $(1+z)^N$.
Therefore, we can conclude that
$\lim_{k\rightarrow\infty}b^{[k]}(z)=b(z)$ for a suitable $b(z)$.
The latter means that, for $k$ large enough,  $\|\bb^{[k]}\|<C$
where $C>0$ is a constant independent of $k$.
Next, we observe that each $\frac {d^\beta \ak(z)}{dz^\beta}$ contains at least $N-\beta$ factors of the form $(1+e^{\lambda_n 2^{-k-1}} z)$ and, since $1-e^{\lambda_n 2^{-k-1}}=O(2^{-k})$, we conclude that
$\frac {d^\beta \ak(-1)}{dz^\beta}=O(2^{-k(N-\beta)})$ as $k\to \infty$. Thus the second part of the claim is proven.
Since the Wronskian matrix ${\cal W}_{\Phi_N}(x)$ is invertible for all $x$ in a neighborhood of zero,
the first part of the claim is obtained from Theorem \ref{Y-THM-10} with $N=1$, so completing the proof.
\endproof

\begin{corollary}
Let $\{\NS,\ k\geq 0\}$ be a $\Phi_N$-generating non-stationary subdivision scheme which
reproduces at least one exponential polynomial in $\Phi_N$ with  ${\cal W}_{\Phi_N}(x)$ invertible for all $x$ in a neighborhood of zero.
If $\{\NS,\ k\geq 0\}$  is asymptotically similar to a convergent stationary scheme, then it satisfies approximate sum rules of order $N$.
\end{corollary}

\begin{remark}\label{rem:notAE}
Note that, the assumption that ${\cal W}_{\Phi_N}(x)$ is invertible for all $x$ in a neighborhood of zero, is needed to prove only
the first of conditions \eqref{SUM-R}.
We additionally emphasize that the assumptions in Theorem \ref{teo:new} do not guarantee asymptotical equivalence between
$\{\NS,\ k\ge 0\}$ and $\{S_\ba\}$. For example, the $k$-level subdivision mask
$$
\begin{array}{l}
\ba^{[k]}=\Big \{ 0, \cdots, 0, \, \frac{1}{(r_k^{-1}+r_k) (r_k^{-\frac12}+r_k^{\frac12})}-\frac{1}{k}, \,
\frac{(r_k^{-\frac12}+r_k^{\frac12})^2-1}{(r_k^{-1}+r_k) (r_k^{-\frac12}+r_k^{\frac12})}-\frac{r_k}{k}, \,\\
\hspace{0.8cm} \frac{(r_k^{-\frac12}+r_k^{\frac12})^2-1}{(r_k^{-1}+r_k) (r_k^{-\frac12}+r_k^{\frac12})}+\frac{1}{k \, r_k^2}, \,
 \frac{1}{(r_k^{-1}+r_k) (r_k^{-\frac12}+r_k^{\frac12})}+\frac{1}{k \, r_k}, \, 0, \cdots, 0 \Big \}, \qquad r_k= e^{\lambda 2^{-k-1}},
\end{array}
$$
is such that $\lim_{k\rightarrow \infty} \| \ba^{[k]} - \ba \|=0$ with $\ba$ denoting the mask of the quadratic B-spline scheme having symbol $a(z)=\frac{1}{4 z^2}(1+z)^3$. The associated symbols satisfy $a^{[k]}(-r_k)=a^{[k]}(-r_k^{-1})=0$ as well as $a^{[k]}(r_k^{-1})=2r_k^{\frac12}$ for all $k \geq 0$, so that $e^{\pm \lambda x}$ are generated whereas only $e^{\lambda x}$ is reproduced with respect to the dual parametrization ${\mathbf t}^{[k]}=\{t_i^{[k]}=\frac{i-\frac{1}{2}}{2^k}, \ i\in \ZZ\}$.
Thus, in view of Theorem \ref{teo:new}, conditions
$| \ak(1) -2|=O(2^{-k})$,
$\left | \ak(-1)\right |
      =O(2^{-2k})$,
$\left | \frac{d \ak(-1)}{dz} \right |
      =O(2^{-k})$, are all satisfied for
$k\to \infty$.
However, since
$\sum_{k=0}^{\infty} \left \vert \frac{1}{(r_k^{-1}+r_k) (r_k^{-\frac12}+r_k^{\frac12})}-\frac{1}{k}-\frac14 \right \vert$
is not convergent, $\{\NS,\ k\geq 0\}$ and $\{S_\ba\}$ are not asymptotically equivalent.
\end{remark}

\smallskip \noindent
We conclude this section by recalling that, in stationary subdivision,
sum rules of order $N$ are known to be
necessary conditions for $C^{N-1}$-continuity
\cite{Cabrelli, ChenJia, DL2002,  Han2010, JetterPlonka}.
As a consequence of the results in \cite{CharinaContiGuglielmiProtasov2014}, in the non-stationary setting,
approximate sum rules of order $N$ and asymptotical similarity to a stationary $C^{N-1}$ subdivision scheme provide sufficient conditions for $C^{N-1}$ regularity of non-stationary subdivision schemes.
More precisely, in view of \cite[Corollary 1]{CharinaContiGuglielmiProtasov2014}, we can state the following.

\begin{corollary}\label{coro:new}
Let $\{\NS,\ k\geq 0\}$ be a $\Phi_N$-generating non-stationary subdivision scheme which
reproduces at least one exponential polynomial in $\Phi_N$. Moreover, let $\{\NS,\ k\geq 0\}$ be asymptotically
similar to a $C^\ell$-convergent subdivision scheme, $\ell\in \NN_0$ and assume that ${\cal W}_{\Phi_N}(x)$ is invertible for all $x$ in a neighborhood of zero.
Then, setting $\rho:=\min\{\ell,N-1\}$, it follows that $\{\NS,\ k\geq 0\}$ is
at least $C^{\rho}$-convergent.
\end{corollary}

\section{Asymptotic behavior of basic limit functions}\label{sec:asymptotic_blf}

It is well-known that a convergent non-stationary subdivision scheme
$\{S_{\ba^{[k]}}, k \geq 0 \}$ defines a family of basic limit functions.
For $\bdelta:=\{\delta_{i,0},\ i\in\ZZ\}$, they are
\begin{equation}\label{phim_def}
\phi_{m}:= \lim_{\ell \rightarrow \infty} S_{\ba^{[m+\ell]}} \cdots S_{\ba^{[m]}} \bdelta, \qquad m \geq 0.
\end{equation}
The goal of this section is to weaken the assumptions used
in \cite[Lemma 15]{DynLevin95} to prove that, as $m\to \infty$,
the sequence of basic limit functions $\{\phi_m,\ m\geq 0\}$
of the non-stationary subdivision scheme $\{S_{\ba^{[k]}}, k \geq 0 \}$
converges uniformly to the basic limit function $\phi$ of a $C^0$
stationary scheme $\{S_{\ba}\}$.  Precisely, in Theorem \ref{teo:conv_phim}
we prove that $\lim_{m \rightarrow \infty} \|\phi_{m} - \phi\|=0$
without requiring that $\{S_{\ba^{[k]}}, k \geq 0 \}$ is
asymptotically equivalent to $\{S_{\ba}\}$.
Differently, we just assume that $\{S_{\ba^{[k]}}, k \geq 0 \}$
is a non-stationary subdivision scheme reproducing one exponential
polynomial in $\Phi_N$, say $e^{\lambda x},\ \lambda\in \CC$,
and it is asymptotically similar to $\{S_{\ba}\}$.
We remark that, the assumption that $\{S_{\ba^{[k]}}, k \geq 0 \}$
reproduces $e^{\lambda x},\ \lambda\in \CC$, is not restrictive.
In fact, if $\{S_{\ba^{[k]}}, k \geq 0 \}$ reproduces a higher order
exponential polynomial $x^r e^{\lambda x},\ r\in \NN$,
then it obviously reproduces $e^{\lambda x}$ as well. Moreover,
if $\{S_{\ba^{[k]}}, k \geq 0 \}$ reproduces constants, with the
choice $\lambda=0$ we can also recover this case.
For the proof of Theorem \ref{teo:conv_phim} we also recall that, accordingly to
\cite[Theorem 1]{ContiRomani11}, a non-stationary subdivision scheme
$\{S_{\ba^{[k]}}, k \geq 0 \}$ reproduces $e^{\lambda x}$, $\lambda \in \CC$
with respect to the parametrization
${\mathbf t}^{[k]}=\{t_i^{[k]}=\frac{i+p}{2^k}, \ i\in \ZZ\}$,
$p \in \{0, - \frac{1}{2}\}$, if its symbol $a^{[k]}(z)$ is such that
$a^{[k]}(-r_k^{-1})=0$ and $a^{[k]}(r_k^{-1})=2r_k^{-p}$, for all $k \geq 0,$
with $r_k= e^{\lambda 2^{-k-1}}$.\\
Before proving Theorem \ref{teo:conv_phim}, we additionally need the
following auxiliary results.

\begin{proposition}\label{prop:nuovaH}
For $p\in \{0,\ -\frac12\}$, let $\{S_{\bh_p^{[k]}}, k \geq 0 \}$ be
the non-stationary subdivision scheme with $k$-level symbol
\begin{equation}\label{SYM-HK}
h_p^{[k]}(z)
:= \frac{r_k^{-p}}{(r^{-1}_k+r_k) z} \, (1+r^{-1}_kz)(1+r_kz), \qquad r_k= e^{\lambda 2^{-k-1}}, \ \ \lambda \in \CC.
\end{equation}
Then 
\begin{itemize}
\item[(a)] for $p\in \{0,\ -\frac12\}$, $\{S_{\bh_p^{[k]}}, k \geq 0 \}$
is $C^0$ and stable (in the sense of \cite{DynLevin95}).
Moreover, the scheme is interpolatory
when $p=0$ while approximating when $p=-\frac12$.
\item[(b)] for $p\in \{0,\ -\frac12\}$, the basic limit functions
\begin{equation}\label{FT-HM}
H_{m,p} :=\lim_{\ell \rightarrow \infty} S_{\bh_p^{[m+\ell]}}
\dots S_{\bh_p^{[m]}} \bdelta, \quad m\ge 0
\end{equation}
of $\{S_{\bh_p^{[k]}}, k \geq 0 \}$
are such that $\lim_{m\rightarrow \infty}H_{m,p}=H$, where $H$ is the linear B-spline supported on $[-1,1]$.
Moreover, $H_{m,0}(0)=1$ while $H_{m,-\frac12}(0)\to 1$ as $m\to \infty$.
\end{itemize}
\end{proposition}

\proof
$a)$ The observation that $|1-r_k| \leq C 2^{-k}$ for some constant $C>0$, implies that
the non-stationary subdivision scheme $\{S_{\bh_p^{[k]}}, k \geq 0 \}$ is asymptotically equivalent to
the linear B-spline scheme, which is $C^0$ and stable (in the sense of \cite{DynLevin95}).
Accordingly, the scheme
$\{S_{\bh_p^{[k]}}, k \geq 0 \}$ is also $C^0$ and stable.
Moreover, by definition of $h_p^{[k]}(z)$ in \eqref{SYM-HK},
$h_p^{[k]}(-r_k^{-1}) = 0$
and $h_p^{[k]}(r_k^{-1}) = 2 r_k^{-p}$, $p\in \{0,\ -\frac12\}$.
Since by definition \eqref{SYM-HK} we can immediately write that $h_p^{[k]}(-z)+h_p^{[k]}(z)=2r_k^{-p}$,
the claim follows.\\
$b)$ Since the basic limit function of the linear B-spline scheme is the degree-1 B-spline $H$ supported on $[-1,1]$,
in view of \cite[Lemma 15]{DynLevin95} we have that $H_{m,p}\to H$ uniformly as $m\to \infty$.
Moreover,
since the mask of the scheme $\{S_{\bh_p^{[k]}}, k \geq 0 \}$
has the same support as the one of the linear B-spline scheme,
we have that for all $m\ge 0$ the basic limit function
$H_{m,p}$ of $\{S_{\bh_p^{[k]}}, k \geq 0 \}$ is supported on $[-1,1]$
and is such that
$H_{m,p}(0) = (e^{\lambda 2^{-m}})^{-p}$. Hence, we can conclude that,
$H_{m,0}(0)=1$ implying that
$H_{m,-\frac 12}(0) \to 1$ as $m\to \infty$.
\endproof

\smallskip \noindent
We continue by providing an additional intermediate result that
will be exploited in the proof of Theorem \ref{teo:conv_phim}.
This result is in Theorem \ref{teo:main}. It extends \cite[Theorem 11]{ContiDynManniMazure14}, where the assumption that $\{\NS, k\geq 0\}$ reproduces constants, is replaced by
the reproduction of at least one exponential polynomial
$e^{\lambda x}$, $\lambda\in \CC$.
For the proof of Theorem \ref{teo:main}, we need to replace the classical notion of backward
difference operator
$(\Delta^{k}\bff^{[k]})_i := \rf_i^{[k]}-\rf_{i-1}^{[k]}$ with
$$(\Delta^{k}_{\lambda} \bff^{[k]})_i := \rf_i^{[k]}-r_{k-1}\rf_{i-1}^{[k]},$$
where $r_{k-1}:=e^{\lambda2^{-k}}$, $\lambda\in \CC$, $k\geq 0$; see also
\cite{UB} for the details of this difference operator.
Obviously, $\Delta^{k}_{\lambda}$ reduces to $\Delta^{k}$ if $\lambda=0$.
The following lemma is useful for our further analysis.\\

Note that, from now on, $C$ will be used to denote any generic positive constant.\\

\begin{lemma}\label{LEM-10}
Let $\{S_{\ba^{[k]}}, k \geq 0 \}$ be a non-stationary subdivision scheme
reproducing $e^{\lambda x},\ \lambda\in \CC$ with respect to the shift parameter $p \in \{0,-\frac{1}{2}\}$, and let
$\{S_{\ba^{[k]}}, k \geq 0 \}$ be asymptotically similar to a convergent stationary scheme $\{S_{\ba}\}$.
For all $m \geq 0$, let
\begin{equation}\label{FM-NS}
\bff_m^{[k+1]}:=
S_{\ba^{[m+k]}}\cdots S_{\ba^{[m]}}\bff^{[0]}.
\end{equation}
Then, there exist a constant $C>0$, $\mu\in (0,1)$ and $K$ large enough
such that
\begin{equation}\label{eq:deltanorm}
\|\Delta_\lambda^{[m+k+1]} \bff_m^{[k+1]}\|\le C \mu^{k} , \quad \forall k \geq K.
\end{equation}
\end{lemma}
\proof
The proof of this result is very similar to  the proof of
\cite[Theorem 3]{ContiDynManniMazure14}. However, there are some
crucial different points that we need to  put in evidence.
From the generation properties of $\{S_{\ba^{[k]}}, k \geq 0 \}$
we know that $a^{[k]}(z)$ factorizes as $a^{[k]}(z)=(1+r_kz) b^{[k]}(z)$
for some Laurent polynomial $b^{[k]}(z)$.
In view of \eqref{F-NS} and \eqref{FM-NS}, we have
$f_m^{[k+1]}(z)=a^{[k+m]}(z) f_m^{[k]}(z^2)$ with $f_m^{[k]}(z):=\sum_{i\in \ZZ}(\bff_m^{[k]})_i z^i\,$.
Multiplying both sides of this equation   by $1-r_{k+m}z$
and using the relation $r^2_{k+m} = r_{k+m-1}$, we arrive at
\begin{align*}
\Delta^{k+m+1}_{\lambda} f_m^{[k+1]}(z)
&=(1-r_{k+m} z)f_m^{[k+1]}(z)
 =(1-r^2_{k+m} z^2)b^{[k+m]}(z)f_m^{[k]}(z^2)\\
& =(1-r_{k+m-1} z^2)b^{[k+m]}(z)f_m^{[k]}(z^2)
  = b^{[k+m]}(z) \Delta^{k+m}_{\lambda} f_m^{[k]}(z^2).
\end{align*}
The latter means that, under the assumption that the  non-stationary
subdivision scheme $\{S_{\ba^{[k]}}, k \geq 0 \}$ generates
$e^{\lambda x}$, there exists a non-stationary subdivision scheme
$\{S_{\bb^{[k]}}, k \geq 0 \}$  such that for the sequence
$\{\bff_m^{[k+1]}:=S_{\ba^{[k+m]}} \bff_m^{[k]},\ k\ge 0\}$,
it is verified that
\begin{equation}\label{DEL-FK}
\Delta^{k+m+1}_{\lambda} \bff_m^{[k+1]}
= S_{\bb^{[k+m]}} \Delta^{k+m}_{\lambda} \bff_m^{[k]},\quad  k\ge 0
\end{equation}
where $\bff_m^{[0]}:=\bff^{[0]}$.
We next show that there exist $\mu\in (0,1)$ and  $K$ large enough
such that, for all $k \ge K$,
$$
\|\Delta^{[k+m+1]}_\lambda\bff_m^{[k+1]}\|\le C \mu^k,
$$
with a constant $C>0$. \\
From the assumption that $\{S_{\ba^{[k]}}, k \geq 0 \}$ reproduces
$e^{\lambda x}$, we know that, in addition to
$a^{[k]}(z)=(1+r_kz) b^{[k]}(z)$, $a^{[k]}(z)$ verifies the condition
$a^{[k]}(r_k^{-1})=2r_k^{-p}$, $p \in \{0,-\frac{1}{2}\}$. Therefore, we can conclude that $b^{[k]}(r^{-1}_k)=r^{-p}_k$.
Moreover, since $\{S_{\ba^{[k]}}, k \geq 0 \}$ is asymptotically similar
to a convergent scheme $\{S_{\ba}\}$, we have that $a(z)=(1+z) b(z)$
with $\lim_{k \rightarrow \infty}\| \bb^{[k]}-\bb\|=0$,
and there exists a positive integer $L$ such that $\|S^L_\bb\|<1$.
Writing the symbol $a^{[k]}(z)$ as $a^{[k]}(z)=\sum_{j=-M}^M \ra_j^{[k]}z^j$,
for the Laurent polynomial $b^{[k]}(z)=\sum_{j=-M}^{M-1} \rb_j^{[k]} z^j$
we have that its coefficients $\rb_j^{[k]}$, $j=-M, \cdots, M-1$ satisfy
\begin{equation}\label{eq:boundmasks}
\rb_{-M}^{[k]}=\ra_{-M}^{[k]},\quad \rb_{j}^{[k]}=\sum_{i=0}^{j-M+1} (-1)^{i} (r^{-1}_k)^{i+1}\rb_{j+1-i}^{[k]},\quad j=-M+1,\dots,M-1.
\end{equation}
Due to $r^{-1}_k\le \max(1, e^{-\lambda})$,
we have
\begin{equation}\label{eq:boundmasks2}
\|\bb^{[k]}\|\le  2MC\|\ba^{[k]}\|\
\end{equation}
with a suitable constant $C>0$.
Since $\lim_{k \rightarrow \infty}\| \bb^{[k]}-\bb\|=0$,
for all $\epsilon>0$, there exists a positive integer $K_\epsilon$ such that
for $k> K_\epsilon$
\begin{equation}\label{eq:S_b}
\|S_{\bb^{[k]}}\|\le \|S_{\bb^{[k]}}-S_{\bb}\|+
\|S_\bb\|<\|S_\bb\|+\epsilon.
\end{equation}
Since for all $k>K_\epsilon$
$$
 \|S_{\bb^{[k+L-1]}}\cdots S_{\bb^{[k]}}-S^L_{\bb}\|\le \sum_{\ell=1}^L \|S_{\bb^{[k+L-1]}}\|\cdots \|S_{\bb^{[k+L-\ell+ 1]}}\|
 \|S_{\bb^{[k+L-\ell]}}- S_{\bb}\| \|S^{L-\ell}_{\bb}\|\,,
$$
due to \eqref{eq:S_b} we get
$$
\|S_{\bb^{[k+L-1]}}\cdots S_{\bb^{[k]}}-S^L_{\bb}\|
\le L \epsilon \left(\|S_{\bb}\|+\epsilon \right)^{\ell -1}\|S_{\bb}\|^{L-\ell}
< C\epsilon,
$$
for a suitable constant $C>0$.
Therefore,
\begin{equation}\label{SB-1}
\|S_{\bb^{[k+L-1]}}\cdots S_{\bb^{[k]}}\|\le \|S^L_{\bb}\|+ C\epsilon,
\qquad \forall k> K_\epsilon.
\end{equation}
Further, since $\|S^L_{\bb}\|<1$, we especially choose a number
$\tilde \epsilon>0$ such that
$${\tilde \mu}:=\|S^L_{\bb}\| +C\tilde \epsilon<1.
$$
Letting $K_{\tilde \epsilon}$ be the corresponding number
satisfying \eqref{eq:S_b},
put $K =\max (K_{\tilde \epsilon}, m)$.
Then,  from \eqref{DEL-FK}, we have that
for an arbitrary integer $n >0$,
$$
 \Delta_\lambda^{[K+m+n+1]} \bff_m^{[K+n+1]}
=S_{\bb^{[K+m+n]}}\cdots S_{\bb^{[K+1]}}
\left(S_{\bb^{[K]}}\cdots S_{\bb^{[m]}}\right)
\Delta_\lambda^{[m]} \bff^{[0]}\,.
$$
Hence, applying (\ref{SB-1}) and using the value of $\tilde \mu$
$$\|\Delta_\lambda^{[K+m+n+1]} \bff_m^{[K+n+1]}\|
\le {\widetilde C} {\tilde \mu}^{\lfloor \frac{n+m} L\rfloor}
\le {\widetilde C} {\tilde \mu}^{\frac {n+m} L-1}
$$
with
$$
{\widetilde C}=\max_{r=0,\cdots,L-1}
\Big(\prod_{j=m}^{K+r}\|S_{\bb^{[j]}}\|\Big)\|
\Delta^{[m]}_\lambda\bff^{[0]}\|.
$$
Moreover, introducing the notation
${C}:={\widetilde C} \,{\tilde \mu}^{\frac {m-K}L-1}$ and
$\mu :={\tilde \mu}^{\frac  1 L}$,
we arrive at the bound
$$\|\Delta_\lambda^{[K+m+n+1]} \bff_m^{[K+n+1]}\|
\le {C} {\mu}^{K+n}.$$
As a  conclusion, there exists a constant $C>0$ such that
$\|\Delta_\lambda^{[m+k+1]} \bff_m^{[k+1]}\|\le C \mu^{k}$ for
$ \mu \in (0,1)$.
\endproof

\begin{theorem}\label{teo:main}
Let $\{S_{\ba^{[k]}}, k \geq 0 \}$ be a non-stationary subdivision scheme
reproducing $e^{\lambda x},\ \lambda\in \CC$ with respect to the shift parameter $p \in \{0,-\frac{1}{2}\}$, and let
$\{S_{\ba^{[k]}}, k \geq 0 \}$ be asymptotically similar to a convergent stationary scheme $\{S_{\ba}\}$.
Then, for all $m \geq 0$, there exist a constant $C>0$, $\mu\in (0,1)$ and $K$ large enough
such that
\begin{equation}\label{eq:Nira}
\left \|
\lim_{\ell \to \infty} \bff_m^{[\ell]}-F_{m,p}^{k} \right \|
\leq C\mu^k \qquad \forall \ k\geq K,
\end{equation}
where $\bff^{[k]}_m$ is in \eqref{FM-NS} and
$F_{m,p}^k$ is the function defined by
\begin{equation}\label{Fmp:def1}
F_{m,p}^{k}:=\sum_{i\in \ZZ}
(\bff_m^{[k]})_i \, H_{m+k,p}(2^{k} \cdot -i),
\end{equation}
with $H_{m+k,p}$  in \eqref{FT-HM}.
\end{theorem}

\proof
Let $\bd^{[k]}=\ba^{[k]}-\bh_p^{[k]}$ with $\bh_p^{[k]}$ the subdivision mask associated with the symbol in \eqref{SYM-HK}.
Since $d^{[k]}(\pm r_k^{-1})=0$ we can obviously write it as
\begin{equation}\label{def:eq:factorD}
d^{[k]}(z)=(1-r_k^2z^2)e^{[k]}(z),
\end{equation}
for a suitable Laurent polynomial $e^{[k]}(z)$.
With the same arguments as in \eqref{eq:boundmasks}
and \eqref{eq:boundmasks2} and using the fact that $\lim_{k\rightarrow\infty} \| \bh_p^{[k]}-\bh\|=0$ in view of Proposition \ref{prop:nuovaH}-(b), we obtain
$$
\|\be^{[k]}\|\le  2M\|\bd^{[k]}\|\le  2M\left(\|\ba^{[k]}\|+\|\bh\|\right)=2M\left(\|\ba^{[k]}\|+1\right)\,.
$$
Moreover, since $\lim_{k\rightarrow\infty}\|\ba^{[k]}-\ba\|=0$, the sequence
$\{\ba^{[k]}, \, k\ge 0\}$ is uniformly bounded independent of $k$, so that
we conclude that
\begin{equation}\label{bounde}
\| \be^{[k]} \| \leq C,
\end{equation}
with $C$ a positive constant.

\smallskip \noindent
The rest of the proof mimics \cite[Theorem 11]{ContiDynManniMazure14}
where $\Delta_\lambda$ replaces $\Delta $ and $H_{m,p}$ replaces
the hat function $H$. By using standard arguments
(see, e.g., \cite[Theorem 4.11]{DL2002}), to show \eqref{eq:Nira}
it suffices to show that, for the sequence of approximating
functions $\{F^{k}_{m,p}, k\geq 0\}$ defined in \eqref{Fmp:def1},
there exists $K>0$ such that for all $k\geq K$,
$\|F_{m,p}^{k+1}-F_{m,p}^{k}\| \leq C\mu^k$ with $\mu\in (0,1)$ and $C$ a positive constant.
For this purpose we write $F_{m,p}^{k}$ as
$$
F_{m,p}^{k}=\sum_{i\in \ZZ}
\left(S_{\bh_p^{[k+m]}}\bff_m^{[k]}\right)_iH_{m+k+1,p}(2^{k+1} \cdot-i)\,.
$$
Hence, it is clear from the relation
$\bff_m^{[k+1]}= S_{\ba^{[k+m]}}\bff_m^{[k]}$
and  the refinability  of $F_{m,p}^{k+1}$ that
$$
F_{m,p}^{k+1}-F_{m,p}^{k} =\sum_{i\in \ZZ}
\left(\left(S_{\ba^{[k+m]}}-S_{\bh_p^{[k+m]}}\right)\bff_m^{[k]}\right)_i
H_{m+k+1,p}(2^{k+1} \cdot-i)\,.
$$
Next, we define the sequence $\bg_m^{[k+1]}$ as
\begin{equation}\label{star}
\bg_m^{[k+1]}:= \left(S_{\ba^{[k+m]}}-S_{\bh_p^{[k+m]}}\right)\bff_m^{[k]}
=S_{\bd^{[k+m]}}\bff_m^{[k]}\,.
\end{equation}
By the expression $g_m^{[k+1]}(z):=\sum_{i\in \ZZ}(\bg_m^{[k+1]})_iz^i$,
we write \eqref{star} as
  $  g_m^{[k+1]}(z)=d^{[k+m]}(z)f_m^{[k]}(z^2)$,
where $f_m^{[k]}(z):=\sum_{i\in \ZZ}(\bff_m^{[k]})_iz^i$.
Then, using \eqref{def:eq:factorD} we have
$$
g_m^{[k+1]}(z)=(1-r_{k+m}^2z^2)e^{[k+m]}(z)f_m^{[k]}(z^2)
=e^{[k+m]}(z)\sum_{i\in \ZZ}
\left((\bff_m^{[k]})_i-r_{k+m-1} (\bff_m^{[k]})_{i-1}\right)z^{2i}\,,
$$
or equivalently
$$
(\bg_m^{[k+1]})_i=\sum_{j\in \ZZ}\re^{[k+m]}_{i-2j}
\left((\bff_m^{[k]})_j-r_{k+m-1} (\bff_m^{[k]})_{j-1}\right)\,,
$$
implying that
$$
|(\bg_m^{[k+1]})_i|
\le \sum_{j\in \ZZ}|\re^{[k+m]}_{i-2j}|\,
   \left| (\bff_m^{[k]})_j-r_{k+m-1} (\bff_m^{[k]})_{j-1} \right|
\le 2M \max_{j}|\re^{[k+m]}_{j}|\,
   \max_{j} \left| \left(\Delta^{[k+m]}_{\lambda}
   \bff_m^{[k]}\right)_{j} \right|\,.
$$
Hence, recalling \eqref{eq:deltanorm} and \eqref{bounde}, there exists $K>0$
such that for all $k \geq K$ we have
$\|\bg_m^{[k+1]}\|\le C\mu^k$ with $\mu \in (0,1)$ and $C$ a positive constant, resulting in that
for the sequence of functions $\{F^{k+1}_{m,p}, k\geq 0\}$
$\|F_{m,p}^{k+1}-F_{m,p}^{k}\| \leq C\mu^k$ for all $k \geq K$ with $\mu\in (0,1)$ and $C$ a positive constant.
Thus we can obtain that for any positive integer $\ell$,
\begin{equation}\label{eq:intermedio}
\| F_{m,p}^{k+\ell}- F_{m,p}^k \|
\leq \sum_{j=0}^{\ell -1} \| F_{m,p}^{k+j+1} -F_{m,p}^{k+j} \|
\leq \sum_{j=0}^{\ell -1} C\mu^{k+j} \leq C \frac{\mu^k}{1-\mu},
\end{equation}
from which it follows that $\{F_{m,p}^{k+1}, k \geq 0\}$ is a Cauchy sequence in the
$L_\infty$-norm and therefore converges uniformly to a continuous limit.
To show that such limit is exactly
$\lim_{\ell \rightarrow \infty}\bff_m^{[\ell]}$
we use standard arguments (see \cite[Lemma 14]{DynLevin95}):
since $\{H_{m+k,p}, \ k \geq 0\}$ is a stable sequence of continuous,
compactly supported functions which approximate partition of unity uniformly,
the uniform convergence of the sequence
$\{F^{k+1}_{m,p}, \ k \geq 0\}$ implies
the uniform convergence of the subdivision scheme
$\{S_{\ba^{[k]}}, k \geq 0 \}$ and also
$$
\lim_{\ell \rightarrow \infty}F_{m,p}^{k+\ell}
=\lim_{\ell \rightarrow \infty} \sum_{i\in \ZZ}
  \left(\bff_m^{[k+\ell ]}\right)_i
  H_{m+k+\ell,p}(2^{k+\ell } \cdot -i)
=\lim_{\ell \rightarrow \infty}
  \bff^{[\ell ]}_m\,,
$$
so concluding the proof.
\endproof

\medskip \noindent
We are finally ready to prove the main result of this section.

\begin{theorem}\label{teo:conv_phim}
Let $\{S_{\ba^{[k]}}, k \geq 0 \}$ be a non-stationary
subdivision scheme reproducing $e^{\lambda x},\ \lambda\in \CC$ with respect to the shift parameter $p\in\{0,-\frac12\}$. Let $\{\phi_m,\ m\geq 0\}$ be the associated sequence of basic limit functions.
Assume further that $\{S_{\ba^{[k]}}, k \geq 0 \}$ is asymptotically similar to a convergent stationary scheme $\{S_{\ba}\}$ with stable basic limit function $\phi$ having H\"{o}lder continuity $\alpha \in (0,1)$.
Then  $\lim_{m \rightarrow \infty} \|\phi_{m} - \phi\|=0$.
\end{theorem}

\proof
Let $k\in \NN_0$ be a non-negative integer,
$\bdelta_m^{[k ]}:=
  S_{\ba^{[m+k-1]}} \cdots S_{\ba^{[m]}} \bdelta$ and $\bphi_m^{[k ]}:=\{\phi_m( i2^{-k}), \, i \in \ZZ \}$.
Then
\begin{equation}\label{eq:product_Soperators}
S_{\ba^{[m+k]}} \cdots S_{\ba^{[m]}} \bdelta - S_{\ba}^{k+1} \bdelta
=\displaystyle{\sum_{j=0}^{k-1}}\, S_{\ba}^j
    \left(S_{\ba^{[m+k-j]}}-S_{\ba} \right)
   \bdelta^{[k-j]}_m
   +  S_{\ba}^k (S_{\ba^{[m]}}-S_{\ba})\bdelta.
\end{equation}
The last term on the right hand side of the above equation can be estimated as
$$
\| S_{\ba}^k  (S_{\ba^{[m]}}-S_{\ba} ) \bdelta\|
\leq \|S_\ba^k\|\, \|\ba^{[m]}-\ba\|.
$$
Due to the convergence assumption of $\{S_\ba\}$,  there exists a constant $C>0$
such that for all $k$, $\|S^k_\ba\|<C$ (see \cite[Section 2]{DynLevin95}).
Hence, from the asymptotical similarity of $\{S_{\ba^{[k]}}, k \geq 0 \}$ and $\{S_{\ba}\}$,
it is immediate that
$$\lim_{m\rightarrow \infty}
\|S^k_\ba\left(S_{\ba^{[m]}}-S_{\ba} \right) \bdelta\|=0.$$
Now, to estimate the other summands in  (\ref{eq:product_Soperators}) note that
\ban\label{AK-A}
((S_{\ba^{[m+k-j]}}-S_{\ba} ) \bdelta_m^{[k-j]})_{i}
  = \sum_{\ell \in \ZZ}
   (\ra^{[m+k-j]}_{i-2\ell}-{\ra}_{i-2\ell} )
   (\bdelta_m^{[k-j]})_{\ell}, \quad i\in \ZZ,\quad j=0,\cdots,k-1.
\ean
To estimate the expression in \eqref{AK-A} we approximate
$\bdelta_m^{[k-j]}$ by using $\bphi_m^{[k-j]}$, \ie \, 
the values of the basic limit functions
$\phi_m$ on the grids $2^{-(k-j)}\ZZ$.
For each $\ell\in\ZZ$ we write
\ban\label{K-DEL}
(\bdelta_m^{[k-j]})_{\ell}
  = \big ((\bdelta_m^{[k-j]})_{\ell} -(\bphi_m^{[k-j]})_{\ell} \big)
   + \big((\bphi_m^{[k-j]})_{\ell} -(\bphi_m^{[k-j]})_i \big)
   + (\bphi_m^{[k-j]})_i.
\ean
In view of (\ref{AK-A}),
 since the masks $\ba^{[k]}$ and $\ba$
have the same finite support,
we need to consider $\ell $ only around $i$, say $|i-\ell|\leq 2M$
for some $M>0$.
Now, to estimate the first term in the right hand side of (\ref{K-DEL}),
since 
$a^{[k]}(r_k^{-1})=2r_k^{-p}$ with $r_k= e^{\lambda 2^{-k-1}}$, we
consider the sequence of functions $\{{\cal F}_{m,p}^{k}, \, k\geq 0\}$ defined by
\begin{equation}\label{FT-FK2}
{\cal F}^k_{m,p}:=\sum_{i\in \ZZ}
(\bdelta_m^{[k]})_i
H_{m+k,p}(2^{k} \cdot -i),
\end{equation}
where $H_{m+k,p}$ is the basic limit function of the non-stationary scheme
with $s$-level symbol
$h_p^{[s]}(z),\ s\geq m+k$
in \eqref{SYM-HK}.
By  Proposition \ref{prop:nuovaH}-(b) we get, for a suitable constant $C$,
$$\big |(\bdelta_m^{[k-j]})_\ell - {\cal F}_{m,p}^{k-j}(\ell 2^{-(k-j)})\big |
= \big |(1-e^{-p \lambda 2^{-m-(k-j)}}) (\bdelta_m^{[k-j]})_\ell \big |
\leq  C  2^{-m-(k-j)} \|\bdelta_m^{[k-j]} \|. $$
This estimate  and   Theorem \ref{teo:main} yield
the bound
\ban\label{K-CONV}
\big|(\bdelta_m^{[k-j]})_{\ell} -\phi_m(\ell2^{-(k-j)})\big|
\leq C( 2^{-(k-j)} +  \mu^{k-j}),\quad \mu\in (0,1),
\ean
for some other constant $C>0$.
Moreover, to estimate $(\bphi_m^{[k-j]})_{\ell} -(\bphi_m^{[k-j]})_i$
with $|i-\ell|\leq 2M$ for some $M>0$,
we exploit the H\"{o}lder continuity of $\phi_m$ which follows by  
\cite[Theorem 3] {CharinaContiGuglielmiProtasov2014} in view of the fact that 
approximate sum rules of order one are satisfied. 
Hence we get the bound
\ban\label{H-CONT}
\big |\phi_m(\ell2^{-(k-j)}) -\phi_m(i2^{-(k-j)})\big|
\leq 2M 2^{-\alpha (k-j)},\quad \alpha\in (0,1).
\ean
Setting $\theta_j:=\big( \mu^{j-1}+2^{-\alpha j}\big)$ we thus 
clearly have that $\{\theta_j, \, j\in \NN\}$ is an absolutely
summable sequence.
Then,
combining \eqref{K-DEL}, \eqref{K-CONV} and \eqref{H-CONT} with
\eqref{AK-A}
and applying the fact that $\|\bphi_m^{[k-j]}\| \leq \|\phi_m\|$,
we have
\begin{align}\label{SS-10}
\begin{split}
\Big | ((S_{\ba^{[m+k-j]}}-S_{\ba} ) \bdelta_m^{[k-j]})_{i} \Big |
& \leq C\theta_{k-j} \sum_{\ell \in \ZZ}\Big |
   (\ra^{[m+k-j]}_{i-2\ell}-{\ra}_{i-2\ell} ) \Big |
   +\Big | (\bphi_m^{[k-j]})_i \sum_{\ell \in \ZZ}
   (\ra^{[m+k-j]}_{i-2\ell}-{\ra}_{i-2\ell} ) \Big | \\
& \leq C  \|\ba^{[m+k-j]} -{\ba} \|
\theta_{k-j}  +
\|\phi_m\| \,
    \Big | \sum_{\ell \in \ZZ} \big (
    \ra^{[m+k-j]}_{i-2\ell}-{\ra}_{i- 2\ell}\big ) \Big |
\end{split}
\end{align}
for some constant  $C>0$.
Now, to bound the term $\big | \sum_{\ell \in \ZZ} \big (
    \ra^{[m+k-j]}_{i-2\ell}-{\ra}_{i- 2\ell}\big ) \big |$
we proceed as follows. Recalling the assumption that $\{S_{\ba^{[k]}}, k \geq 0 \}$ reproduces $e^{\lambda x}$ with $\lambda \in \CC$, we can apply Theorem \ref{Y-THM-10} in the case $N=1$, to write that
\begin{equation}\label{eq:Y-THM-10_N1}
\Big \vert \ak(1) -2 \Big \vert = \Big \vert \sum_{\ell \in \ZZ} \ra^{[k]}_{2\ell} + \sum_{\ell \in \ZZ} \ra^{[k]}_{2\ell+1} -2 \Big \vert = O(2^{-k}), \quad
k\to \infty.
\end{equation}
Then, using the fact that $\{S_{\ba}\}$ is a convergent stationary scheme, we also have that
$\sum_{\ell \in \ZZ} \ra_{2\ell}+ \sum_{\ell \in \ZZ} \ra_{2\ell+1}=2$.
As a consequence, from  \eqref{eq:Y-THM-10_N1} we easily obtain that
$$
\Big \vert  \sum_{\ell \in \ZZ} \Big( \ra^{[k]}_{2\ell} - \ra_{2\ell} \Big ) + \sum_{\ell \in \ZZ} \Big (  \ra^{[k]}_{2\ell+1} - \ra_{2\ell+1} \Big) \Big \vert = O(2^{-k}),
\quad  k\to \infty.  $$
Similarly, due to the fact that $\ak (-1) = O(2^{-k})$ as $k\to \infty$
(see Theorem \ref{teo:main}),
we get
$$ \Big \vert \sum_{\ell\in \ZZ}\Big ( \ra^{[k]}_{2\ell} - \ra_{2\ell} \Big )
-\sum_{\ell \in \ZZ}
\Big ( \ra^{[k]}_{2\ell+1} - \ra_{2\ell+1} \Big)
\Big \vert = O(2^{-k}),
\quad  k\to \infty,
$$
implying that, for all $i \in \ZZ$,
$$
\Big \vert  \sum_{\ell \in \ZZ} \Big( \ra^{[k]}_{i-2\ell} - \ra_{i-2\ell} \Big ) \Big \vert = O(2^{-k}),
\quad  k\to \infty.
$$
Hence, from the latter equation the bound $\Big|\sum_{\ell \in \ZZ} ( \ra^{[m+k-j]}_{i-2\ell}-{\ra}_{i-2\ell} ) \Big|
\leq C 2^{-(m+k-j)}$ is obtained straightforwardly.
Moreover, since we also have the bound $\|S^j_\ba\|<C$ for all $j$, for the summation in (\ref{eq:product_Soperators}) we are finally able to write
from \eqref{SS-10} that
\begin{align}
\begin{split}
\displaystyle{\sum_{j=0}^{k-1}}
\big|( S^j_\ba(S_{\ba^{[m+k-j]}}-S_{\ba} )
     \bdelta_m^{[k-j]} )_i\big|
\leq C \displaystyle{\sum_{j=0}^{k-1}
\Big( } \|\ba^{[m+j+1]}-{\ba }\|
\theta_{j} + 2^{-(m+j+1)} \Big )
\end{split}
\end{align}
with a suitable constant $C>0$.
Now we let $m\rightarrow \infty$. From the asymptotical similarity
of $\{\ba^{[k]},\ k\geq  0\}$ and $\{\ba\}$, we know that for all
$\epsilon>0$ there exists $\bar m$ such that for $m+j>\bar m $ we have
$\|\ba^{[m+j]}-{\ba }\|\leq \epsilon$. Therefore we can conclude that $\lim_{m\rightarrow \infty}
\|S_{\ba^{[m+k]}} \cdots S_{\ba^{[m]}} \bdelta - S_{\ba}^{k+1}\bdelta\| =0$.
\endproof

\section{ Approximation order}\label{Sec:approx}
In this section we estimate the approximation order of
a $\Phi_N$-reproducing non-stationary subdivision scheme $\{ \NS,\ k\geq 0 \}$
in case the initial data are sampled from a function in the Sobolev space $W^n_\infty(\RR)$, $n\in \NN$.
The latter is defined to be the set of all functions
$f$ in $L_\infty(\RR)$ that have derivative $\frac{d^\ell f}{d u^\ell} \in L_\infty(\RR)$
for all $0\leq \ell \leq n$. We recall that, for any $f\in W^n_\infty(\RR)$, the associated norm
is defined by
\ban\label{F-NORM}
 \|f\|_{n,\infty}: = \sum_{\ell=0}^n \left\|\frac{d^\ell f}{d u^\ell} \right\|_{L_\infty(\RR)}.
\ean
Let $f\in W^n_\infty(\RR)$ and let 
the initial data be of the form
$\bff^{[m]}:=\{\rf_i^{[m]}=f(2^{-m}i),\ i\in\ZZ\}$ for some $m \geq 0$. In the following theorem we estimate the convergence order of the
error  $\|f -g_{\bff^{[m]}}\|_{L_\infty(\RR)}$ as $m\to \infty$, where $g_{\bff^{[m]}}$ is the limit of the subdivision scheme obtained from the initial data $\bff^{[m]}$. We emphasize that
Theorem \ref{TH-AO} extends
the result in \cite[Theorem 2.4]{JKLY}, where the more restrictive assumption of asymptotical equivalence rather than asymptotical similarity is assumed.

\begin{theorem} \label{TH-AO}
Assume that the non-stationary scheme $\{ \NS,\ k\geq 0\}$ is $\Phi_N$-reproducing and is asymptotically similar to a convergent stationary scheme $\{S_\ba\}$. Assume further that the initial data are of the form
$\bff^{[m]}:=\{\rf_i^{[m]}=f(2^{-m}i),\ i\in\ZZ\}$ for some fixed $m\geq 0$ and for some function $f\in W^{\gamma}_\infty(\RR)$ where $\gamma\in \NN$, $\gamma \leq N$.
If the Wronskian matrix ${\cal W}_{{\Phi}_\gamma}(0)$ of ${\Phi}_\gamma \subseteq \Phi_N$ is invertible,
then
\baa
\|g_{\bff^{[m]}}  -f \|_{L_\infty(\RR)} \leq C_f  2^{- \gamma m }\, \quad m\ge 0
\eaa
with a constant $C_f>0$ depending only on  $f$.
\end{theorem}
\proof
For the Sobolev exponent $\gamma$, let ${\Phi}_\gamma := \{\vp_0,\ldots, \vp_{\gamma-1}\}$, $\gamma \leq N$ and let $x$ be a fixed point in $\RR$.
Our proof employs an auxiliary
function $\psi_x$, which depends on $x$ and is defined by
$$\psi(u) := \psi_x(u) :=\sum_{n=0}^{\gamma-1} d_n \varphi_n (u-x),\quad u\in \RR$$
where the entries of the coefficient vector ${\bf d} :=(d_n,\ n =0,\cdots, \gamma-1)$
are obtained by solving the linear system
\begin{equation}\label{HY-04}
\frac{d^r\psi(x)}{du^r}= \frac{d^r f(x)}{du^r}, \quad r =0,\dots, \gamma-1.
\end{equation}
In this way, $\psi$ equals $f$ at $x$, but not on $\RR$.
Equations in \eqref{HY-04} can be
equivalently written in the matrix form
\ban\label{WP-20}
{\cal W}_{{\Phi}_\gamma} (0)\cdot{\bf d}^T = {\bf f}^T,
\ean
with ${\bf f} :=(\frac{d^r f(x)}{du^r}, \ r =0,\cdots, \gamma-1)$. Due to
the invertibility of the Wronskian matrix
${\cal W}_{{\Phi}_\gamma}(0)$,
the uniqueness of the solution of this linear system is guaranteed.
Moreover, since the function $\psi$ belongs to $\Phi_N$ and
the non-stationary scheme $\{\NS, k \geq 0\}$ is $\Phi_N$-reproducing, we can write
$$\psi(u)=\sum_{i \in \ZZ} \psi(2^{-m}i)\phi_m(2^m u -i),\quad u\in \RR,$$
with $\{\phi_m,\ m\geq 0\}$ denoting the basic limit functions
of $\{\NS, k \geq 0\}$  defined in \eqref{phim_def}.
By assumption, $\rf_i^{[m]}=f(2^{-m}i)$, $i\in\ZZ$ and, in view of the linearity of the subdivision operators, we can write
\begin{equation}\label{eq:gf0}
g_{\bff^{[m]}}(u)=\sum_{i\in\ZZ} \rf^{[m]}_i \phi_m(2^m u-i),\quad u\in \RR.
\end{equation}
Thus, using the expression of $g_{\bff^{[m]}}$ in \eqref{eq:gf0} and observing that $f(x)=\psi(x)$
due to the construction of $\psi$ in (\ref{HY-04}), we can rewrite $f(x) -g_{\bff^{[m]}}(x)$ in the following way:
$$
f(x) - g_{\bff^{[m]}}(x)
 = \psi(x)-\sum_{i \in  \ZZ } f(2^{-m}i) \phi_m(2^m x-i)\\
 = \sum_{i \in  \ZZ} \big (\psi(2^{-m}i)-f(2^{-m}i)\big ) \phi_m(2^m x-i)
$$
Now, let $T^\gamma_g(u):= \sum_{\ell =0}^{\gamma -1}(u-x)^\ell \frac{d^\ell g(x)}{du^\ell}/\ell !$ be the degree-$(\gamma-1)$ Taylor polynomial
of a function $g\in C^{\gamma-1}(\RR)$  around $x$.
Then, consider the Taylor expansions $T^\gamma_\psi$ and
$T^\gamma_f$ of the functions $\psi$ and $f$, respectively.
Due to the condition in (\ref{HY-04}),
it is obvious that $ T^\gamma_{\psi} (2^{-m}i)  =T^\gamma_{f} (2^{-m}i).$
Hence, applying the remainder form of the Taylor expansion, we get
\begin{align*}
|f(x) -g_{\bff^{[m]}}(x)|
   &\leq \frac{2^{-m \gamma}}{\gamma!} \sum_{i\in\ZZ } \left|   (i-2^m x)^\gamma \, \frac{d^\gamma\, (f-\psi)(\xi_i)}{d u^\gamma}\phi_m(2^m x-i) \right|,
\end{align*}
for some $\xi_i$ between $x$ and $i2^{-m}$.
By (\ref{HY-04})  and \eqref{WP-20}, we are able to write
$|\frac{d^\gamma \psi(\xi_i)}{d u^\gamma}|\leq C_\gamma \|f\|_{\gamma,\infty}$
for some positive constant $C_\gamma$
independent of $x$ and $i2^{-m}$, 
with $\|\cdot \|_{\gamma,\infty}$ in \eqref{F-NORM}.
Thus,  it is immediate that
\begin{align}\label{HY-80}
|f(x) -g_{\bff^{[m]}}(x)|
   &\leq C_\gamma \frac{2^{-m \gamma}}{\gamma!} \|f \|\low{\gamma,\infty}
        \sum_{i\in\ZZ }|\phi_m(2^m x-i)
        (i-2^m x)^\gamma|.
\end{align}
Finally, by Theorem \ref{teo:conv_phim}, $\phi_m$ is uniformly bounded independent of $m$.
Moreover, since $\phi_m$ is compactly supported,
$\#\{i\in \ZZ:\phi_m(2^m x -i)\not =0\}\leq C$ for any fixed $x$ and $m$.
Therefore the claim follows from (\ref{HY-80}).
\endproof

\section{Conclusions}\label{conclusions}
In this paper we have shown that for primal or dual non-stationary subdivision schemes (essentially the ones of interest in applications) the reproduction of $N$ exponential polynomials implies approximate sum rules of order $N$. This mimics the stationary case where the reproduction of $N$ polynomials implies sum rules of order $N$.
Furthermore, also in analogy to the stationary case where generation of $N$ polynomials implies sum rules of order $N$, we have shown that generation of $N$ exponential polynomials implies approximate sum rules of order $N$ if asymptotical similarity to a convergent stationary scheme is assumed together with reproduction of a single exponential polynomial.\\
We additionally considered the non-stationary counterpart of the well-known result asserting that the reproduction  of an $N$-dimensional space of polynomials is sufficient for the subdivision scheme to have approximation order $N$.
In particular, for a non-stationary subdivision scheme, the reproduction
of an $N$-dimensional space of exponential polynomials, jointly with asymptotical similarity, has been shown to imply approximation order $N$.
The proof of the latter required us to show also the uniform convergence of the sequence of basic limit functions of a non-stationary scheme,  reproducing one exponential polynomial, to the basic limit function of the asymptotically similar stationary scheme.\\
We finally remark that, since asymptotical similarity is needed to get several of our results, we cannot claim that approximate sum rules are the complete satisfactory notion to replace sum rules when moving from the stationary to the non-stationary setting. However, this notion is definitely very helpful in almost all practical cases.

\medskip
\section*{Acknowledgements}
Support from the Italian GNCS-INdAM is gratefully acknowledged. Lucia Romani acknowledges the support received from Ministero dell'Istruzione, dell'Universit\`a e della Ricerca - Progetti di Ricerca di Interesse Nazionale 2012 (MIUR-PRIN 2012 - grant 2012MTE38N).
Jungho Yoon acknowledges the support received from the National Research Foundation of Korea through
the grants NRF-2015-R1A5A1009350 (Science Research Center Program) and NRF-2015-R1D1A1A09057553.

\end{document}